\documentclass[11pt]{article}
\usepackage{amsmath}
\usepackage{amssymb}
\usepackage{amscd}
\usepackage[dvips]{graphics}
\usepackage{amsfonts}

\def\reel{\hbox{{\rm R}\kern-1em\hbox{{\rm I} }}}
\def\relatif{\ \hbox{{\rm Z}\kern-.4em\hbox{\rm Z}}}
\def\nat{\hbox{{\rm N}\kern-1em\hbox{{\rm I} } }}
\def\comp{\hbox{{\rm C}\kern-.55em\hbox{{\rm I} } }}
\def\smallcomp{\hbox{\fiverm C}\kern-.35em{\hbox{\fiverm I}}}

\def\fudge{\mathchoice{}{}{\mkern.5mu}{\mkern.8mu}}
\def\bbc#1#2{{\rm \mkern#2mu\vbar\mkern-#2mu#1}}
\def\bbb#1{{\rm I\mkern-3.5mu #1}} \def\bba#1#2{{\rm
#1\mkern-#2mu\fudge
#1}}
\def\bb#1{{\count4=`#1 \advance\count4by-64 \ifcase\count4\or\bba
A{11.5}\or
\bbb B\or\bbc C{5}\or\bbb D\or\bbb E\or\bbb F \or\bbc G{5}\or\bbb H\or
\bbb I\or\bbc J{3}\or\bbb K\or\bbb L \or\bbb M\or\bbb N\or\bbc O{5}
\or
\bbb P\or\bbc Q{5}\orrrr\b
bb R\or\bbc S{4.2}\or\bba T{10.5}\or\bbc U{5}\or    \bba V{12}\or\bba
W{16.5}\or\bba X{11}\or\bba Y{11.7}\or\bba Z{7.5}\fi}}
\def\rat{\hbox{{\rm Q}\kern-.70em\hbox{{\rm I} } }}

\newcommand{\refm [1]} {(\ref{#1})}

\newcounter{theorem}

\newtheorem{theorem}{Theorem}[section]
\newtheorem{lemma}[theorem]{Lemma}

\newtheorem{thm}[theorem]{Theorem}

\newtheorem{proposition}[theorem]{Proposition}

\newcommand{\be}{\begin{equation}}
\newcommand{\eu}{\end{equation}}
\newcommand{\ber}{\begin{eqnarray}}
\newcommand{\ena}{\end{eqnarray}}
\newcommand{\nin}{\noindent}
\newcommand{\non}{\nonumber}

\def\qed{\hfill \vrule height1.3ex width1.2ex depth-0.1ex}
\def\bbb#1{{\rm I\mkern-3.5mu #1}} \def\bba#1#2{{\rm
#1\mkern-#2mu\fudge
#1}}

\newcommand{\xx}{\bf x}
\newcommand{\sr}{\scriptscriptstyle}
\newcommand{\ds}{\displaystyle}

\newcommand{\dN}{\delta}

\newcommand{\nun}{\nu_{\sr{N}}}

\newcommand{\la}{\label}

\title{Reversible coagulation-fragmentation processes and random combinatorial structures: asymptotics for
the number of groups}

\author{{\bf Michael M. Erlihson}
\thanks{E-mail: maerlich@tx.technion.ac.il}\\
Department of Mathematics, Technion-Israel Institute of Technology,\\
Haifa, 32000, Israel.\\
\quad
{\bf Boris L. Granovsky}
\thanks{E-mail: mar18aa@techunix.technion.ac.il} \\
Department of Mathematics, Technion-Israel Institute of Technology,\\
Haifa, 32000, Israel.}
\begin{document}
\maketitle \vskip 5cm \nin American Mathematical Society 2000
subject classifications.

\nin Primary-60J27; secondary-60K35, 82C22, 82C26.

\nin Keywords and phrases: Coagulation-Fragmentation processes,
Random combinatorial structures, Local and central limit theorems,
Distributions on the set of partitions.
\newpage

\begin{abstract}
  The equilibrium distribution of a reversible coagulation-fragmentation process (CFP) and the joint distribution of
  components of a random combinatorial structure(RCS) are given by the same probability measure on the set of
  partitions. We establish a central limit theorem for the number of groups (=components) in the
  case $a(k)=qk^{p-1},\ k\geq1,\ q,p>0$, where $a(k),\ k\geq1$ is the parameter function that induces the invariant measure.
  The result obtained is compared with the ones for logarithmic RCS's and for RCS's, corresponding to the case $p<0$.
\end{abstract}

 \section{ Summary.}
  \nin Our  main result is a central limit theorem (Theorem $4.6$) for the number of groups at steady state
  for a class of reversible CFP's and for the corresponding class of RCS's.

  \nin In Section 2, we provide a definition of a reversible k-CFP admitting interactions of up to $k$ groups,
   as a generalization of the standard 2-CFP. The steady state of the processes considered is fully defined by a parameter
   function $a\ge 0$ on the set of integers. It was observed by F. Kelly (\cite{kel}, p. 183) that for
   all $2\le k\le N$ the k-CFP's have the same invariant measure on the set of partitions of a given
   integer $N(=$ the number of particles).

  \nin Section 3 explains the idea of A. Khintchine's probabilistic method for derivation of asymptotic formulae.
  In the spirit of the method, we construct a representation of the probability function of the number of groups via the
  probability function of the sum of i.i.d. random variables.

  \nin In section 4 we study the case when the parameter function $a$ is of the
  form: $a(k)=qk^{p-1},\ k\geq1,\\ q,p>0.$  We prove a local and a central normal  limit theorems for
  the number of groups at equilibrium, as $N\to\infty$. To achieve this, we employ a new (for this field) tool:
  the Poisson summation formula. We conclude the section by providing some intuition for the main result of the paper.

  \nin In Section 5 we recall that the invariant measure of a reversible CFP can be viewed as a joint distribution
  of components of a RCS, known as an assembly. Comparing our main result with the known ones for RCS's, we identify $p=0$
  as a point of phase transition of the invariant measure, as $N\to\infty$. We also provide a few examples of RCS's that
  conform to the case $p>0$.

 \section{CFP's with multiple interactions: Definition. }

 \nin  Following \cite{kel}, \cite{dgg}, we treat a CFP  as a continuous-time Markov chain on the finite set
  $\Omega_N=\{\eta\}$ of all partitions  $\eta=(n_1,\ldots,n_N)$ of a given integer $N$:\be \sum_{j=1}^{N} j n_j =N ,
  \quad n_j\ge 0, \quad j=1,\ldots,N.\la{**}\end{equation} Here $N$ codes the total population of indistinguishable
  particles partitioned into groups(=clusters) of different sizes. A  group  of size $j\ge 2$ may split into
  a number, say $s,\ 2\le s\le j$, of groups of sizes $j_1, \ldots,j_s: j_1+\ldots+j_s=j$ and, conversely, the above $ s$
  groups may coagulate into one large group of size $j$. We will call these $s$-interactions(=s-transitions), $J_s$ -fragmentation
  and $J_s$- coagulation respectively, where $J_s=(j_1,\ldots,j_s).$  A
  stochastic process that admits interactions of up to $k\leq N$ groups will be denoted k-CFP. Note that both types of the
  interactions conserve the total number of particles.

  \nin We now  provide a formal definition of a  k-CFP that naturally extends the definition of  the
  standard 2-CFP (see \cite{dgg}, \cite{kel}). A k-CFP is given by the rates of infinitesimal transitions
  that are assumed to depend only on the sizes of interacting groups. For  given $N$ and $s,\;2\le s\le k$,
  let $\psi_s,\phi_s\ge 0 $ be a pair of functions defined on the same set $$\mathbb{J}_{s}=\{J_s=(j_1,\ldots,j_s): j_l>0,
  \quad l=1,\ldots ,s,\quad j_1+\ldots j_s\le N\}$$ \nin of s-tuples of positive integers, depicting the sizes of
  interacting groups. The functions $\psi_s,\phi_s$ are the rates of $J_s$-coagulation and $J_s$-fragmentation respectively.
  Both functions are assumed to be invariant w.r.t. all $s!$ permutations of $j_1,\ldots, j_s.$
  To complete the definition of a  k-CFP, it is left to determine the total rates of all possible s-transitions from one
  partition $\eta\in\Omega_N$ to another. Assume that the given $J_s=(j_1,\ldots,j_s)\in \mathbb{J}_{s}$ and
  $\eta=(n_1,\ldots,n_N)\in\Omega_N$ are such that $n_l\geqslant m_l\ge1,\  l\in J_s$, where $m_l$ counts the number
  of components in $J_s$ that are equal to $l$. In other words, we assume, that a $J_s$-coagulation of some
  groups in the partition $\eta$ is possible. Clearly, a given $J_s$-coagulation of any groups in $\eta$ transforms $\eta$
  into the same partition that will be denoted $\eta^{(J_s)}\in\Omega_N.$ Similarly, the result of $J_s$-fragmentation of
  any groups in $\eta\in\Omega_N$ will be denoted $\eta_{(J_s)}\in\Omega_N$. By a simple combinatorial calculation,
  the total rate $\Psi_s(J_s;\eta)$ of all possible $J_s$- coagulations at the partition $\eta$ is

 \be \Psi_s(J_s;\eta)= \bigg(\prod_{j_l}
 \frac{n_{j_l}!}{(n_{j_{l}}-m_{l})!}\bigg)\psi_s(J_s), \quad
 J_{s}\in \mathbb{J}_{s}, \la{01}
 \end{equation}

 \nin where the product is taken over all distinct components $j_l$ of $J_s$. By the same logic we define
 the total rate $\Phi(J_s;\eta)$ of all possible $J_s$-fragmentations  at the partition $\eta$ to be equal to
 \be
  \Phi_s(J_s;\eta)=n_{\vert J_s\vert}\phi_s(J_s), \quad J_s\in
  \mathbb{J}_{s}, \la{02}
 \end{equation}
 \nin where $\vert J_s\vert=j_1+\ldots+ j_s .$ Now we see that a k-CFP is fully defined by the $(k-1)$ pairs of
 functions $\psi_s,\phi_s,\quad s=2,\ldots,k.$

  \nin  In this paper we will be concerned only with reversible CFP's. Define the ratio of s-interactions

\be
q_s(J_s)=
  \begin{cases}
   \frac{\psi_s(J_s)}{\phi_s(J_s)}, &\text{if}\quad
    J_s:\psi_s(J_s)\phi_s(J_s)>0, \quad J_s\in \mathbb{J}_{s},\\
    0, & \text{otherwise}.
  \end{cases}
  \la{04}
\end{equation}
 \nin
  A natural extension of Theorem 1 in \cite{dgg} gives  the following characterization of reversible k-CFP's with
  positive transition rates.

\begin{proposition}

 \nin Let
 \be
  \psi_s(J_s)\phi_s(J_s)>0,\quad J_s\in \mathbb{J}_{s},\
  s=2,\ldots,k.
  \la{--}
 \end{equation}
 \nin Then the corresponding k-CFP is reversible iff the ratios $q_s,\ s=2,\ldots,k$ are of the form:
 \be
 q_s(J_s)=\frac{a(\vert J_s\vert)}{\prod_{l=1}^sa(j_l)},\quad J_s=(j_1,\ldots,j_s)\in
 \mathbb{J}_{s},\quad s=2,\ldots, k,
 \la{03}
\end{equation}
 \nin where $a>0$ is a function on the set of positive integers.
 \end{proposition}

 \nin The proof is deferred until after Proposition 2.2.

 \vskip .5cm \nin {\bf A historical remark:} It was already noted in \cite{dgg}, that the characterization of reversible
 2-CFP's was motivated by the following two completely independent lines of research: the seminal paper of F. Spitzer
 \cite{sp} on nearest-particle systems and F. Kelly's and P. Whittle's works in the 1970-s on networks and clustering
 process in polymerization (see \cite{kel},
 \cite{wh2}).\qquad\qquad\qquad\qquad\qquad\qquad\qquad$\blacksquare$
 \vskip .5cm

 \nin We will call $a$ the parameter function of a reversible CFP and we will write $a_j=a(j),\ j\geq1$. The following result
 gives the explicit form of the steady state of the processes.

\vskip .5cm
 \begin{proposition} (\cite{dgg},\cite{kel},Ch.8)

\nin  For a given $N$,  all k-CFP's, $k=2,\ldots ,N$ satisfying
\refm[03] for some parameter function $a$ have the same invariant
measure $\mu_N$ on $\Omega_N:$

\begin{center}
\be
 \mu_N(\eta)=({c_N})^{-1}\:\frac{a_1^{n_1}
 a_2^{n_2}\ldots a_N^{n_N}} {{n_1}!{n_2}!\ldots{n_N}!},\qquad
 \eta=(n_1,\ldots,n_N)\in \Omega_N,
\la{05}
\end{equation}
\end{center}
 where $a_j>0 ,\;j=1,2,\ldots,N$ and $c_N=c_N(a_1,\ldots,a_N)$ is the partition function of the measure $\mu_N:$
 \be
  c_0=1,\quad
  c_N=\sum_{\eta\in\Omega_N}\frac{a_1^{n_1} a_2^{n_2}\ldots  a_N^{n_N}}{{n_1}!{n_2}!\ldots{n_N}!}\,,\quad
  \quad N\geq1.
  \la{08}
 \end{equation}
 \end{proposition}
 \nin {\bf Proof:} $\quad$ We have to show that the measure $\mu_N$ given by $\refm[05]$ satisfies the detailed
 balance condition:
 \be
 \la{401}
 \mu_N(\eta)V(\eta,\xi)=\mu_N(\xi)V(\xi,\eta),\qquad \eta,\:\xi\in\Omega_N,
 \end{equation}
 where $V(\eta,\xi)$ is the total rate of the infinitesimal (in time) transition from
 $\eta$ to $\xi$. If $V(\eta,\xi)=0$, then \refm[401] is trivially
 true. If now $\xi=\eta^{(J_s)}$ for some $J_s\in\mathbb{J}_s,\:2\leqslant s\leqslant k$, then we see from
 \refm[05] and \refm[01] that
 \be
  \frac{\mu_N\big(\eta^{(J_s)}\big)}{\mu_N(\eta)}=\frac{\Psi_s(J_s;\eta)}
  {\Phi_s(J_s;\eta^{(J_s)})}.
 \end{equation}
 \nin In a similar manner one can verify \refm[401] in the case
 $\xi=\eta_{(J_s)}$ for some $J_s\in\mathbb{J}_s$.\qquad\qquad$\blacksquare$\\
 \nin Now we are in position to give the
 \\\nin{\bf Proof} of Proposition 2.1:\quad If a k-CFP, $k\geq2$ is
 reversible, then it follows from $\cite{dgg}$, Theorem 1, that \refm[03] should hold for
 s=2. The latter fact implies that the unique invariant measure of all k-CFP-s, $k=2,\ldots,N$ is given by \refm[05].
 Consequently, \refm[03] should hold for $3\leq s\leq k$, by Proposition 2.2. The same reasoning proves the
 converse part of the claim. \qquad$\blacksquare$
 \\ \nin The preceding discussion shows that the steady state of a reversible k-CFP is uniquely determined  by a parameter
function $a$.

\section{Khintchine's type representation for the probability function
of the number of groups.}

\nin Our objective will be the study of the asymptotic behavior,
as $N\to \infty,$ of the number of groups $\nu_N$ at equilibrium
given by the measure $\mu_N$. It follows from \refm[05] and
\refm[08] that
 \begin{equation}
 \mathbb{P}\:(\nu_{\scriptscriptstyle{N}}=n)={(c_N)}^{-1}
 \bigg(\sum_{\eta\in\Omega_N}\frac{a_1^{n_1} a_2^{n_2}\ldots
a_N^{n_N}}
 {{n_1}!{n_2}!\ldots{n_N}!}\mathbf{1}_{(\sum_{k=1}^N{n_k}=n)}\bigg),\
\quad n\leq N. \la{09}
 \end{equation}

 \nin As in \cite{frgr1}, \cite{frgr2}, our tool will be the
 probabilistic method by A. Khintchine introduced in the 1950's in
 his seminal book \cite{Kh} . The idea of the method is to
 construct  the representation of the quantity of interest via
 the probability function of a sum of independent integer valued
 random variables depending on a free parameter, and subsequent implementation of a local limit theorem.
 This allows for the derivation of the desired asymptotic formula.

 \nin In number theory, the implementation of Khintchine's method
 was developed by G. Freiman \emph{et al} (for references see
 \cite{post},\cite{frgr2}). In particular, a general scheme of the
 method for asymptotic problems related to partitions was
 outlined by G. Freiman and J. Pitman in \cite{fr}. The method was applied to CFP's for the first time
 in \cite{frgr1}, for derivation of the asymptotic formula for
 the partition function of the measure $\mu_N$ in the case $a_k\sim
 k^{p-1},\ k\to \infty, \ p>0.$ In \cite{frgr2} the method was used for the study
 of the asymptotic behavior of some quantities related to
 clustering of groups at the steady state, when $a_k\sim
 k^{p-1}L(k), \ k\to \infty,\ p>0,$ where $L$ is a slowly varying function.

  \nin Though the implementation of Khintchine's method goes along the standard scheme, the related
 asymptotic analysis varies from problem to problem. In contrast
 to the aforementioned research, the problem treated in the
 present paper requires knowledge of the  second term   in the
 asymptotic expansions considered. In light of this, we employ the Poisson summation formula,
 a new tool for this field.

 \nin We will always assume that the parameter functions $a$ considered
 are positive and s.t. the power series
 \be
  \sum_{k=1}^{\infty}{{a_k}x^k}, \quad x\in \mathbb{C}
  \la{000}
 \end{equation}
 \nin has a finite radius of convergence $R>0$. Since
 the transformation $a_k\Rightarrow h^k a_k, \ h>0,\ k=1,2,\ldots
 ,N$ does not change the measure $\mu_N$ given by \refm[05], we
 assume w.l.g. that $R=1.$
 \\ \nin Now let $\xi_1,\ldots,\xi_n$ be i.i.d. integer valued nonzero
 r.v.'s defined by
  \be
   \mathbb{P}({\xi_1}=l)=\frac{{a_l e^{-{\delta}l}}}{S(\delta)},\quad a_l>0,  \quad l\geq1,
   \la{1}
  \end{equation}
 \nin where $\delta>0$ is a free parameter and
 \be
  S(\delta)=\displaystyle\sum_{k=1}^{\infty}{{a_k}e^{-{\delta}k}}.
  \la{2}
 \end{equation}
 \nin Note that the r.v. $\xi_1$ has finite moments of all orders
 for all $\delta>0,$ since $R=1$.
 \nin  We start with the following representation of
 the probability $\mathbb{P}({\nu_N}=n).$

 \begin{lemma}
 \nin Define $$ T_n=\sum_{k=1}^{n}{\xi_k},\quad n\geq1, $$ where
  $\xi_k,\ k=1,2,\ldots $ are i.i.d. r.v.'s given by \refm[1],\refm[2].
 Then
 \ber
  \mathbb{P}({\nu_{\sr{N}}}=n)={(c_N n!)}^{-1}{S^n(\delta)}{e^{{\delta}N}}
  {\mathbb{P}({T_n}=N)}, \quad \delta>0.
  \la{3}
 \ena
 \end{lemma}
\nin {\bf Proof:} It follows from \refm[1],\refm[2] that
 \be
   \mathbb{P}({T_n}=N)=\frac{\displaystyle\bigg(\sum_{\eta\in\Omega_N}\frac{a_1^{n_1}
  a_2^{n_2}\ldots a_N^{n_N}}{{n_1}!{n_2}!\ldots{n_N}!}\mathbf{1}_{(\sum_{k=1}^N
  {n_k}=n)}\bigg){n!}}{{S^n(\delta)}{e^{{\delta}N}}}.
 \end{equation}
 \nin By \refm[05] and \refm[09] this implies the claim \refm[3].\qed
 \vskip.5cm
 \nin {\bf Remark.}  \refm[3] has a form of a typical
 representation in Khintchine's method. It can be also viewed as a
 version of the representation formula for the total number of
 components in the generalized scheme of allocation (see
 \cite{kol}, Lemma 1.3.3).
 \vskip.5cm
  \nin Our study will be heavily based on the assumption
   \be
    \la{0011}
    \sum_{k=1}^{\infty}{ka_k}=\infty.
   \end{equation}
  Denote
 \be
   M_1=M_{1}(n;\delta)=\mathbf{E}{T_n}=n\mathbf{E}{\xi_1}=nS^{-1}(\dN)
  \sum_{k=1}^{\infty}ka_k{e^{-{\delta}k}}, \quad \dN>0 \la{4a}
 \end{equation}

 \nin and choose the free parameter $\delta$ as a solution of the
 equation

 \be
   M_1= N,\quad n\leq N.
  \la{5}
 \end{equation}
 \nin Such a choice of the free parameter is typical for Khintchine's method (\cite{Kh}, p.110) and is designed to
 make the probability ${\mathbb{P}({T_n}=N)}$ in \refm[3] large, as $n,N\to \infty, \ n\le
 N.$

 \begin{lemma}
  \nin Under condition \refm[0011], the equation \refm[5] has a unique solution $\delta=\delta_{n,N}$ for all  $n\leq
  N.$
 \end{lemma}

 \nin{\bf Proof}: We first show  that $M_1$ is decreasing in $\delta>0:$
 \be
   M_1^\prime(n;\delta)=n\ \frac{-S(\delta)\sum_{k=1}^{\infty}k^2a_ke^{-\delta k}+ \Big(\sum_{k=1}^{\infty}ka_k
   e^{-\delta k}\Big)^2}{S^2(\delta)}<0,\quad \delta>0,
 \end{equation}
  \nin by the Cauchy-Schwarz inequality. Consequently, $M_1(n;0)=\ds{\sup_{\delta>0}{M_1(n;\delta)}:=nM_1^{*}}$
  where $M_1^{*}$ does not depend on $n$. If now $S(0)<\infty$, then \refm[0011] immediately implies
  $M_1^{*}=\infty.$ In the case $S(0)=\infty$, supposing $M_1^{*}<\infty$ leads to the contradiction:
  \ber
   \non 0&>& M_1(n;\delta)-nM_1^{*}=n\frac{\sum_{k=1}^{\infty}{a_ke^{-\delta
   k}(k-M_1^{*})}}{S(\delta)}\\&\geq& nS^{-1}(\delta)\left(A(\delta)+\frac{1}{2}\sum_{k\geq2M_1^{*}}{ka_ke^{-\delta
   k}}\right),\quad\delta>0,
  \ena
  where $A(\delta)$ is bounded for any $\delta\geq0$, while the sum in the brackets tends to $+\infty$ as
  $\delta\to0^{+}$, by \refm[0011]. Hence, $M_1(n;0)=+\infty$. Finally,
  \be
   M_1(n;\infty)=n\ds\lim_{\delta\to\infty}\frac{\sum_{k=1}^{\infty}{ka_ke^{-\delta(k-1)}}}
   {\sum_{k=1}^{\infty}{a_ke^{-\delta(k-1)}}}=n.
  \end{equation}
 \qed
 \section{A central limit theorem for the number of groups.}
 \nin  Our paper is devoted exclusively to the case when the parameter function $a$ has a polynomial rate of
 growth, namely:
 \be
 a_k=qk^{p-1}, \quad q,p>0,\quad k\geq1.
  \la{001}
 \end{equation}

 \nin We first consider the case $q=1.$
 The following lemma which is basic for our subsequent asymptotic
  analysis is a particular case of the Poisson
 summation formula (see \cite{br},\cite{chandr}).

\begin{lemma}(\cite{chandr}, p.82).  Let $p>1$ and $Re(z)>0.$ Then we have

 \be
  \sum_{k=1}^{\infty}{e^{-zk}k^{p-1}}=\Gamma(p)\bigg(\sum_{l=-\infty}^{\infty}({{z+{2\pi}il}})^{-p}\bigg).
  \la{10}
 \end{equation}
 \end{lemma}

 \nin With the help of this remarkable identity we derive the following asymptotic
 formula that holds for $p>0$.

 \begin{lemma}
   If $p>0,\ |z|\rightarrow {0},\ \ Re(z)\rightarrow{0^{+}},$ then
 \be
   \sum_{k=1}^{\infty}{e^{-zk}k^{p-1}}=\Gamma(p)\Big(z^{-p}+A(p)\Big)+O(z),
   \la{11}
 \end{equation}

 \nin where $A(p)$ is a constant which in the case $p>1$ is given
 explicitly by
 \be
   A(p)=2{(2{\pi})}^{-p}\zeta(p)\cos{\frac{{\pi}p}{2}}. \la{12}
 \end{equation}
 \nin (Here $\zeta(p)$ is the Riemann zeta function).
 \end{lemma}
 \nin {\bf Proof}: First consider the case $p>1.$ Let
 \be
   G(z)=\sum_{l=-\infty}^{-1}{(z+2{\pi}il)}^{-p}+\sum_{l=1}^{\infty}
   {(z+2{\pi}il)}^{-p},\quad Re(z)\ge 0, \quad p>1.
   \la{12a}
 \end{equation}
 \nin The two series in the RHS of \refm[12a] converge absolutely
 when $Re(z)>0,$ while, by straightforward calculations,
 $G(0)=A(p),\ G^\prime(0)=-pA(p+1)$ with $A(p)$ given by \refm[12].
 Consequently, by \refm[10] we have for $p>1$
 \ber
  \non\sum_{k=1}^{\infty}{e^{-zk}k^{p-1}}-\Gamma(p)z^{-p}&=&
  \Gamma(p)G(z)=\Gamma(p)\bigg(G(0)+G^\prime(0)z+o(z)\bigg),\\
  & &|z|\rightarrow {0},\ \ Re(z)\rightarrow{0^{+}}.
  \la{12b}
 \ena
 \nin This proves \refm[11] for $p>1.$  So, we write
  for $p>0$

 \be
  \sum_{k=1}^{\infty}{e^{-zk}k^{p}}=\Gamma(p+1)\Big(z^{-p-1}+A(p+1)\Big)+O(z),
  \quad  |z|\rightarrow {0},\ \ Re(z)\rightarrow{0^{+}}. \la{12c}
 \end{equation}
 \nin Next, integrating \refm[12c] w.r.t. $z$ gives \refm[11] with
  a constant $A(p)$ that is not known explicitly.\qed

 \nin Now we are in a position to derive the asymptotic formula for
  the free parameter $\delta$.

 \begin{proposition}
  \nin Assume that $n$ is s.t.\ \  $\alpha:=\frac{N}{n}\to \infty,$ as $ N\to \infty$. Then
 \be
  \delta=p{\alpha}^{-1}\Big({1-A(p){p}^{p}{\alpha}^{-p}}\Big)+o\big({\alpha}^{-p-1}\big),\quad \alpha\to \infty.\la{13}
 \end{equation}
 \end{proposition}
 \nin {\bf Proof}: First observe that in the case of the function $a$ considered, it follows from
 \refm[4a],\refm[5] and the proof of Lemma 3.2, that $\alpha=\frac{N}{n}\to \infty$ implies $\delta_{n,N}\to
 0^+$. Implementing \refm[11]  gives

 \ber
   \mathbf{E}{\xi_1}&=&\frac{\displaystyle\sum_{k=1}^{\infty}k^{p}{e^{-{\dN}k}}}
   {\displaystyle\sum_{k=1}^{\infty}k^{p-1}{e^{-{\dN}k}}}=p\
   \frac{\ds{{\dN}^{-p-1}}
   +A(p+1)+O(\delta)}{\ds{{\dN}^{-p}}+A(p)+O(\delta)}= \non \\
   &=&p{\ds{\dN}^{-1}}{\Big(1-A(p){\dN}^{p}+o({\dN}^{p})\Big)},\quad
  \delta\to 0^{+},\quad p>0.
 \ena
 \nin Now \refm[5] leads to

 \be
  \delta=p\alpha^{-1}{\Big(1-A(p){\dN}^{p}+o({\dN}^{p})\Big)},\quad
  \dN=\dN_{n,N}\to 0^+, \quad \alpha\to \infty.
 \end{equation}

 \nin Iterating this equation w.r.t. $\delta$ gives \refm[13].\qed
 \vskip.5cm
 \nin We will focus now on asymptotics, as $\alpha\to \infty,$ of the probability ${\mathbb{P}({T_n}=N)}$
  under $\dN$ given by \refm[13]. We have
 \be
  \mathbb{P}({T_n}=N)={(2\pi)}^{-1}\bigg(\int_{-\pi}^{\pi}{\varphi(t)e^{-itN}}dt \bigg), \la{14}
 \end{equation}

 \nin where $\varphi$ is the characteristic function of $T_n.$

 \nin The basic idea of the Khintchine's method is that for a wide
  class of models,  choosing  the free parameter from the condition
 \refm[5] guarantees that the main contribution to the integral
  in the RHS of \refm[14] comes from a set which is some
  neighborhood of zero.

 \nin We will demonstrate that this is in effect true in the case
 considered and will prove the normal local limit theorem for the
 sum $T_n,$ as $\alpha\to \infty.$ As a preliminary step, we verify
 the validity of the Lyapunov's condition

 \be
  \lim_{n\to\infty}\frac{M_3}{ M_2^{\frac{3}{2}}}=0, \la{4f}
 \end{equation}
 \nin where $M_2=M_2(n;N),\ M_3=M_3(n;N)$ are correspondingly the
  variance and the third central moment of $T_n$ under $\dN$ given
  by \refm[13].

 \nin We have
 \be
  M_2=Var{T_n}=nVar{\xi_1}=n\Big(S^{-1}(\delta)
  \sum_{k=1}^{\infty}k^2a_k{e^{-{\delta}k}}-\alpha^2\Big), \quad \delta>0
  \la{4}
 \end{equation}
 and
 \be
  M_3=E(T_n-N)^3=nE(\xi_1-\alpha)^3=n(E\xi_1^3-3\alpha E\xi_1^2+2\alpha^3). \la{4b}
 \end{equation}
\nin Applying now \refm[11] and \refm[13] gives the following
asymptotic expressions for the moments considered:
 $$
  M_2\thicksim\ n\bigg(\ p(p+1){\dN}^{-2}-p^2{\dN}^{-2}\bigg)
  =np{\dN}^{-2},$$
 \be
   M_3\thicksim\ n\bigg(\ p(p+1)(p+2){\dN}^{-3}-3p^2(p+1)\delta^{-3}+2p^3\delta^{-3}\bigg)=2p{\dN}^{-3}n,
  \quad \ n,\alpha\to \infty.
  \la{4c}
 \end{equation}

 \nin This proves \refm[4f]. The condition \refm[4f] provides the existence of  $\beta=\beta(n,\delta)>0$ s.t.
 \be
  \beta^2 M_2\to \infty, \quad \beta^3 M_3\to 0, \quad n,\alpha\to\infty. \la{4g}
 \end{equation}
 \nin Explicitly, in view of \refm[4c],
 \be
  \beta=\delta n^{-\frac{1-\epsilon}{2}}, \ 0<\epsilon<1/3
  \la{4n}
 \end{equation}
 \nin satisfies \refm[4g]. As it will be shown below, $[-\beta,\beta]$ is just the required neighborhood of zero.

 \nin
 \begin{lemma}
  (The local limit theorem for $T_n$.)

 \nin Set
 \be
   n=Q_pN^{\frac{p}{p+1}}+sN^{\frac{p}{2p+2}},\quad Q_p=p^{-1}(\Gamma(p+1))^{\frac{1}{p+1}},\quad s\in R,\quad p>0.
   \la{40}
 \end{equation}

 \nin Then
 \be
  \mathbb{P}({T_n}=N)\sim \frac{1}{\sqrt{2\pi Var T_n}}\ ,\quad N\to\infty. \la{17}
 \end{equation}
 \end{lemma}
 \nin {\bf Proof}: We write

 \be
  I=I_1+I_2,
  \la{16}
 \end{equation}

 \nin where $I=\ds{\int_{-\pi}^{\pi}{\varphi(t) e^{-itN}}dt}\ $ and $I_1,I_2$ are integrals of the integrand of $I$
 taken over the sets $[-\beta,\beta ]$ and $[-\pi,-\beta] \bigcup [\beta,\pi]$ respectively.

 \nin {\bf Step 1}. We  find the asymptotics of the integral $I_1,$ when
  $\beta$ is as given by \refm[4n].

 \nin  By the definition of $\alpha,$
 \be
  \varphi(t) e^{-itN}=\varphi_1^n(t),\quad t\in R, \la{18}
 \end{equation}
 \nin where $\varphi_1$ is the characteristic function of the r.v.
 $\xi_1-\alpha.$ Next, denoting by $\varphi_2$ the characteristic function of the r.v.
 $\left(p^{-1/2}\delta\right)(\xi_1-\alpha),$  we get from \refm[4]--\refm[4c]

 \be
  \varphi_2(t)=1- \frac{1}{2}t^2+O(t^3),\quad t\to 0.
  \la{19}
 \end{equation}
 \nin Combining this with the relationship
 \be
  \varphi_2(t)=\varphi_1(p^{-1/2}\delta t),\quad t\in R, \la{20}
 \end{equation}
 \nin  \refm[18] becomes
 $$ \varphi(t) e^{-itN}=\bigg(1- \frac{1}{2}(\sqrt{p}{\dN}^{-1}t)^2+O(\delta^{-3}t^3)\bigg)^n\sim
 $$
 \be
  \exp \bigg(- \frac{1}{2}(\sqrt{np}{\dN}^{-1}t)^2+nO(\delta^{-3}t^3)\bigg), \quad t\delta^{-1}\to 0,\quad
  n\to\infty.
  \la{21}
\end{equation}

 \nin Consequently, by virtue of \refm[4c] and \refm[4g],
 \be
   I_1\sim \sqrt{\frac{2\pi}{VarT_n}}, \quad \quad n,\alpha\to\infty,\quad p>0.
   \la{22}
 \end{equation}
 \nin {\bf Step 2}. We are to show that

 \be
  I_2=o(I_1), \quad  N\to \infty, \quad p>0.
  \la{0101}
 \end{equation}

 \nin We apply \refm[10] with $z=\delta- it , \quad t\in[\beta,\pi]$ to obtain from \refm[12a] the following analog of
 \refm[12b]:
 \be
  \sum_{k=1}^{\infty}{e^{-zk}k^{p-1}}=\Gamma(p)\bigg(z^{-p}+G(-it)+\delta G^\prime(-it)\bigg)+o(\delta), \quad p>1,
  \quad t\in[\beta,\pi], \quad \delta\to 0^{+},
  \la{222}
 \end{equation}
 \nin where $\vert G(it)\vert, \vert G^\prime(it)\vert\le const:=B(p),\quad t\in [\beta,\pi],\quad p>1.$ By the same
 argument as for \refm[12c] the latter yields
 \be
  \vert\varphi _1(t)\vert\le\frac{(\delta^2+\beta^2)^{-\frac{p}{2}}+B(p)+O(\delta)}{\delta^{-p}+A(p)+
  O(\delta)},\quad t\in [\beta,\pi], \quad p>0, \quad \delta\to 0^{+}.
  \la{223}
  \end{equation}

 \nin In the rest of this section it is always assumed that $n$ is specified as in \refm[40]. In view of \refm[40]  we have
 \be
   n\alpha^{-p}=\frac{n^{p+1}}{N^p}=p^{-p}\Gamma(p)+ \epsilon_N,\quad p>0, \quad N\to \infty, \la{4111}
 \end{equation}

 \nin where
 \be
  \epsilon_N=(p+1) Q_p^p sN^{-\frac{p}{2p+2}}+\frac{p(p+1)}{2}Q_p^{p-1}
  s^2N^{-\frac{p}{p+1}}+O(N^{-\frac{3p}{2p+2}}),\quad p>0, \quad N\to \infty. \la{400}
 \end{equation}

 \nin This fact will be repeatedly used in our subsequent asymptotic analysis.  Consequently,
 \be
   n\delta^p= \Gamma(p)+p^p\epsilon_N+O\big(\alpha^{-p}\big), \quad p>0,\quad N\to\infty.
   \la{4200}
 \end{equation}

 \nin Employing \refm[4200] we further obtain from \refm[223] and \refm[4n]
  \be
  \vert\varphi _1(t)\vert^n\le O\bigg(\exp{(-\frac{p}{2}n^\epsilon)}\bigg),\quad t\in [\beta,\pi], \quad
  p>0, \quad n\to \infty,
 \end{equation}

 \nin for any $0<\epsilon<\frac{1}{3}$. This together with \refm[22] proves \refm[0101]. \qed

 \nin To establish our main result, it is left to find the asymptotic formulae for the rest of the factors in \refm[3].
 First we make use of the following result of \cite{frgr1}:

 \be
   c_N\sim\frac{1}{\sqrt{2\pi B_N^2}}\exp\bigg(N\sigma+
   \sum_{k=1}^N k^{p-1} e^{-k\sigma}\bigg),\quad N\to \infty,\quad  p>0,
   \la{2444}
 \end{equation}
 \nin where
 \be
  \sigma=\sigma_N\sim\bigg(\frac{N}{\Gamma(p+1)}\bigg)^{-\frac{1}{p+1}},\quad p>0,\quad N\to \infty
  \la{27}
 \end{equation}
 \nin  is the unique solution of the equation
 \be
  \sum_{k=1}^{N}e^{-k\sigma}k^p=N,\quad p>0,
  \la{25}
 \end{equation}
 \nin while
 \be
  B_{N}^{2}=\sum_{k=1}^{N}e^{-k\sigma}k^{p+1},\quad p>0.
  \la{26}
 \end{equation}

  \nin For our purpose, we need to know the second term in the asymptotic expansion \refm[27] of $\sigma$. In
  what follows we denote by $\epsilon_N$ different quantities tending to zero, as $N\to \infty.$ By our
  asymptotic formula \refm[11],

 \be
   \ds{\sum_{k=1}^{\infty}k^{p}e^{-k\sigma}}=\Gamma(p+1) \bigg(\sigma^{-(p+1)}+A(p+1)\bigg)+O(\sigma), \quad p>0,
   \quad N\to\infty,
   \la{28}
 \end{equation}
 \nin where $A(p)$ is as in Lemma 4.2, while, by the Euler summation formula, we have for any $\gamma> 0,$
 and $\sigma=\sigma_N$ given by \refm[27],

 \be
  \sum_{k=N+1}^{\infty}k^{p}e^{-k\sigma}\sim \int_{N+1}^\infty x^pe^{-\sigma
  x}dx=\sigma^{-(p+1)}\int_{\sigma(N+1)}^{\infty}{x^pe^{-x}dx}=o(N^{-\gamma}),\quad p>0, \quad N\to
  \infty.
  \la{29}
 \end{equation}

 \nin In view of \refm[28],\refm[29], the equation  \refm[25] can be rewritten now as
 \be
  \Gamma(p+1) \bigg(\sigma^{-(p+1)}+A(p+1)\bigg)+O(\sigma)=N.
  \la{299}
 \end{equation}

\nin Consequently, we get
$$
 \sigma=\bigg(\frac{N}{\Gamma(p+1)}-A(p+1)+\epsilon_N\bigg)^{-\frac{1}{p+1}}=
$$
  \be
  \bigg(\Gamma(p+1)\bigg)^{\frac{1}{p+1}}N^{-\frac{1}{p+1}}+\frac{A(p+1)}{p+1}
  \bigg(\frac{\Gamma(p+1)}{N}\bigg)^{\frac{p+2}{p+1}}+o(N^{-\frac{p+2}{p+1}}),
  \quad p>0, \quad N\to\infty. \la{30}
 \end{equation}

\nin This yields
 \be
   N\sigma=\ds{N^{\frac{p}{p+1}}}(\Gamma(p+1))^{\frac{1}{p+1}}+\epsilon_N, \quad p>0,\quad N\to \infty.\la{32A}
 \end{equation}
 \nin Next, we obtain from \refm[11],\refm[26],\refm[29] and \refm[30]
 \be
   B^2_N \sim\left(\frac{N}{\Gamma(p+1)}\right)^{\frac{p+2}{p+1}}\Gamma(p+2),\quad p>0, \quad N\to \infty
   \la{32B}
  \end{equation}
\nin and  $$
 \sum_{k=1}^N k^{p-1} e^{-k\sigma}= \sum_{k=1}^\infty k^{p-1}e^{-k\sigma}+\epsilon_N=
 $$
 \be
   \Gamma(p)\bigg(\frac{N}{\Gamma(p+1)}\bigg)^{\frac{p}{p+1}}+\Gamma(p)A(p)+\epsilon_N,\quad p>0,\quad N\to \infty.
   \la{31}
 \end{equation}
 \nin  Substituting the above  expressions in \refm[2444] gives the desired asymptotic formula for $c_N$:
 \be
   c_N\sim h_1N^{-\frac{p+2}{2(p+1)}}\exp(h_2\ds{N}^{\frac{p}{p+1}}+h_3),\quad p>0, \quad N\to \infty,\la{31c}
 \end{equation}
 \nin where the constants $h_i,\ i=1,2,3$ are  given by $$
  h_1=\frac{\bigg(\Gamma(p+1)\bigg)^{\frac{1}{2(p+1)}}}{\sqrt{2\pi(p+1)}},
 $$
  $$ h_2=(p+1)Q_p, $$
 \be
   h_3=\Gamma(p)A(p), \quad p>0. \la{32}
 \end{equation}

  \nin Next, \refm[4111], \refm[40] and \refm[13] give
 \be
  \exp(\delta N)\sim\exp{\bigg(pn-A(p)\Gamma(p+1)\bigg)}, \quad p>0,\quad N\to
  \infty,
  \la{244}
 \end{equation}

 \nin while \refm[4200] gives
 \ber
   \non S^n(\delta)&=&\bigg(\Gamma(p)\delta^{-p}+A(p)\Gamma(p)+O(\delta)\bigg)^n\sim\bigg(\Gamma(p)\bigg)^n\delta^{-pn}
   \exp\bigg(A(p)\Gamma(p)\bigg),\\ & &\quad p>0,\quad N\to \infty.
   \la{245}
 \ena

 \nin We again apply \refm[4111] to get

 $$ \delta^{-pn}\sim p^{-pn}\alpha^{pn}\exp\bigg(A(p)\Gamma(p+1)\bigg),$$

 \be
  \alpha^{pn} \sim \frac{n^np^{pn}\exp{\bigg(-n\Big(\frac{p^p}{\Gamma(p)}\epsilon_N-\frac{1}{2}\big(\frac{p^{p}}
  {\Gamma(p)}\big)^{2}\epsilon_N^2+O(\epsilon_N^3)\Big)\bigg)}} {\bigg(\Gamma(p)\bigg)^n},\quad p>0,\quad N\to \infty,
  \la{246}
 \end{equation}

 \nin where $\epsilon_N$ is given by \refm[400]. Observing that $n\epsilon_N^3\to 0,\quad N\to \infty,$ we write out
 now the asymptotic expressions for $n\epsilon_N$ and $ n\epsilon_N^2$ to obtain

  $$ S^n(\delta)\sim n^n\exp{\bigg(-\frac{s^2}{2d_p}-s(p+1)N^{\frac{p}{2p+2}}+A(p)\big(\Gamma(p)+\Gamma(p+1)\big)\bigg)}, $$

  \be
   d_p=\frac{Q_p}{p+1},\quad p>0,\quad N\to \infty.
   \la{46}
 \end{equation}

 \nin Finally, substituting  in \refm[3] the preceding asymptotic expansions and employing Stirling's asymptotic formula
 gives the following result.

 \nin
  \begin{thm} (The local limit theorem for $\nu_N$).

  \nin Let $n$  be given as in \refm[40]. Then
  \be
   \mathbb{P}(\nu_{\sr{N}}=n)\sim\frac{1}{\sqrt{2\pi d_p}}N^{-\frac{p}{2p+2}} \exp\bigg(-\frac{s^{2}}{2d_p}\bigg):=f(N;s),
   \quad s\in R,\quad p>0,\quad N\to \infty.
   \la{47}
 \end{equation}
 \end{thm}
 \vskip.5cm
 \nin This leads to our main result that says that almost all the mass of the probability distribution of the r.v.
  $\nu_N$ is concentrated , as $N\to \infty ,$ in a neighborhood of size $O(N^{\frac{p}{2p+2}})$ of the point
  $Q_pN^{\frac{p}{p+1}}$.

 \nin
 \begin{thm}
 (The central limit theorem for $\nu_N$).

 \nin
 \be
   \frac{\nu_{\sr{N}}- Q_p\ds{N^{\frac{p}{p+1}}}}{\sqrt{d_p}N^{\frac{p}{2p+2}}}\Rightarrow N(0,1) , \quad p>0,
   \quad N\to\infty,
   \la{48}
 \end{equation}
 \nin where $\Rightarrow$ denotes weak convergence and $d_p,Q_p$ are  as in \refm[46],\refm[40] respectively.

 \end{thm}
 \nin {\bf Proof}: We provide a sketch of the proof that is done by the implementation of the standard technique
 of passing from the local theorem to the integral theorem (for more details see (\cite{shir}, p.59, \cite{dur}, p.81).

 \nin It follows from \refm[47] that for any $a\le b, \quad a,b\in  R,$
 \be
  \mathbb{P}\bigg(\frac{\nu_{\sr{N}}- Q_p\ds{N^{\frac{p}{p+1}}}}{N^{\frac{p}{2p+2}}}\in [a,b]\bigg)=\sum_{s\in R_{N}} f(N;s),
  \quad p>0,\quad N\to \infty, \la{49}
 \end{equation}

 \nin where $f(N;s)$ is given by \refm[47] and $$ R_N=\{s\in[a,b]:n=Q_p\ds{N^{\frac{p}{p+1}}}+
 sN^{\frac{p}{2p+2}}\in \mathbb{N}\}. $$

 \nin Since $\vert R_N\vert=O( N^{\frac{p}{2p+2}}), $ as $N\to\infty,$ we have
 \be
  \sum_{s\in R_N}{exp\bigg(-\frac{s^2}{2d_p}\bigg)}N^{-\frac{p}{2p+2}}\rightarrow\int_a^b{\exp\bigg({-\frac{x^2}{2d_p}}\bigg)
  dx},\quad p>0, \quad N\to \infty.\quad \qed
  \la{490}
 \end{equation}
 {\bf Remark}
  : Since the transformation  $a_k\Rightarrow h^{k}a_k,\;h>0,\; k\geq1$  of the parameter function $a$ does not change
  the measure $\mu_N$, the results of our paper are true for  $a_k=h^kk^{p-1},\ k\geq1,\ h>0,\
  p>0$.\quad$\blacksquare$

 \nin One more extension of Theorem 4.6 is provided by the following result.
 \begin{thm}
 : If $\widetilde{a_k}=qk^{p-1},\ p>0,\ k\geq1$, where $q>0$ is a constant, then
 \ber
   \frac{\nu_{\sr{N}}-\widetilde{Q_p}\ds{N^{\frac{p}{p+1}}}}{\sqrt{\widetilde{d_p}}N^{\frac{p}{2p+2}}}\Rightarrow N(0,1),
   \quad p>0,\quad  N\to \infty,
   \la{480}
 \ena
   where $\widetilde{Q_p}=q^{\frac{1}{p+1}}Q_p\:$, $\widetilde{d_p}=q^{\frac{1}{p+1}}d_p$.
 \end{thm}
   {\bf Proof}: Denote by $\widetilde{\bullet}\ $ the quantities related to the parameter function $\widetilde{a}$. We see
   from \refm[1],\refm[25] and \refm[30] that the r.v. $\widetilde{\xi_1}$ has the same distribution as $\xi_1$ and that
   \be
     \widetilde{\sigma_N}\sim q^{\frac{1}{p+1}}\sigma_N,\quad
     N\to\infty.
     \la{0001}
    \end{equation}
  Repeating the preceding asymptotic analysis gives the claimed change in the scaling induced by
   $q>0$.\qquad $\blacksquare$ \\\\
  \nin We conclude this section by providing some intuition for the scaling in the central limit Theorem 4.6. For this purpose we employ
  the following result established in \cite{frgr2} for the model in question. Denote by $\bar{q}(\eta),\underline{q}(\eta)$
  the size of the largest (resp. smallest) group in a random partition $\eta\in\Omega_N. $ Then
 \ber
   \la{60a}
   \lim_{N\to\infty}{\mathbb{P}\left(N^{\frac{1}{p+1}-\epsilon}<\bar{q}(\eta)<N^{\frac{1}{p+1}+\epsilon}\right)}=1,\quad
   p>0,
 \ena
  for all $\epsilon>0$, while
  \be \
  \ds{\lim_{N\to \infty} \mathbb{P}\big(\:\underline{q}(\eta)\ge l\:\big)}=
  \begin{cases}
    0, & \text{if}\quad l=N^\beta, \quad 0< \beta\le 1,\\
    \exp\Big(-\ds\sum_{j=1}^{l-1} a_j\Big) ,  & \text{if}
     \quad l\ge 2 \quad \text{is a fixed number}.
   \end{cases}
   \la{61b}
 \end{equation}\nin
  From \refm[60a] and \refm[61b] one may conjecture, that for large $N$, the main "mass",  of size $O(N)$, is partitioned
  into groups(=clusters) of sizes $O\left(N^{\frac{1}{p+1}}\right)$, while the rest of the mass, of size $o(N)$, is  partitioned into
  groups of small sizes. Adopting this conjecture gives the expectation of the number of groups
  $\nu_N$ as $O\left(\frac{N}{N^{\frac{1}{p+1}}}\right)=O\left(N^{\frac{p}{p+1}}\right)$.
   Concerning the asymptotic variance of $\nu_N$, we note that the relationship $Var(\nu_N)=O(\mathbb{E}\nu_N),\ N\to\infty$
  holds also for permutations, the Ewens sampling formula and for some other RCS's (see \cite{bar},  \cite{kol}).

 \nin \section{RCS's: examples and comparison with known results.}
 \setcounter{equation}{0}

 \nin  It was already  observed in  \cite{frgr1} and \cite{frgr2} that particular cases of the expression \refm[05]
 for the equilibrium measure $\mu_N$ conform to joint distributions of components of a variety of RCS's, known as
 assemblies.

 \nin  By a combinatorial structure (CS) of a size $N$ one  means a union of nondecomposable
 components(=components) of different sizes. Formally, such a structure is
  given by the two sets of integers $\{p_N,\ N\geq1\}$ and $\{m_N,\ N\geq1\}$
   that count respectively
 the total number of instances of size
 $N$ and the number of components of size $N.$ An example of a CS is a graph on $N$ vertices treated as a union
 of its connected components. Therefore, an instance of a CS of size $N$ is given by a partition $\eta=(n_1,\ldots, n_N)\in
 \Omega_N$, where $n_k$ is the number of components of size $k$ in the instance. By assuming that for a given $N$ an
 instance is chosen randomly from all $p_N$ instances, one induces a RCS that is
 completely determined by the random component size
 counting process, the latter being  a random vector with values in $\Omega_N.$ With an obvious abuse of
 notation  we denote the random vector by $\eta.$ A remarkable fact in the theory of RCS's is that a great variety of
 them obey the conditioning relation
 \be
   {\cal L}(\eta)={\cal L} (Z_1,\ldots, Z_N\vert \sum_{j=1}^N jZ_j=N),
   \la{90}
 \end{equation}

 \nin where $Z_1,\ldots,Z_N$ are independent integer valued r.v.'s. In particular, a class of CS's known as assemblies is
 characterized by the fact that $Z_k$ are Poisson r.v.'s \  s.t. $Z_k\sim Po(m_k/k!), \ k=1,2,\dots, N.$ Consequently, by a
 straightforward calculation we find from \refm[90] that for assemblies
 \be
  \mathbb{P}(\eta)=\mu_N(\eta), \quad \eta\in \Omega_N,
  \la{91}
 \end{equation}
 \nin with $a_k=\frac{m_k}{k!}$ and $c_k=\frac{p_k}{k!},\quad k=1,2,\ldots,N$. In particular, the case $m_k\sim \theta(k-1)!,\;\theta>0$ conforms to
  logarithmic RCS's that encompass permutations ($m_k=(k-1)!$) and the Ewens sampling formula ($m_k=\theta(k-1)!,\quad\theta>0)$.
 \nin The novel theory of general logarithmic RCS's is presented in \cite{bar}. From the analytical
 point of view, the common feature of logarithmic RCS's is that they do not obey the condition $p>0$ adopted in the
 present paper.
 As a result, our asymptotic analysis  based on the Poisson summation formula is not applicable for logarithmic RCS's.

 \nin The case of permutations has a long history. For this case, the seminal result by V.L. Goncharov (1944)
 and L. Shepp and  S. Lloyd (1966) states the following central limit theorem for $\nu_N$:

 \be
  \frac{\nun- \log N}{\sqrt{\log N}}\Rightarrow N(0,1).\quad
  \la{50}
 \end{equation}
 \nin A version of \refm[50] for general logarithmic RCS's is also known.
 \\ L.  Mutafchiev  (\cite{mut}) proved a local limit theorem  for $\nu_N$ under some assumptions on the asymptotic
  behavior as  $x\rightarrow1$ of the generating function of the sequence $\{c_n\}_{1}^{\infty}$ . It can be verified
  that in the case $p>0$ the assumptions $(2.7)$ and $(2.8)$ in \cite{mut} do not hold. However, Mutafchiev
  (\cite{mut},p.425) conjectured, that a result similar to his Theorem 2.4 holds for a wider class of RCS's. Our Theorem 4.6
  confirms this conjecture. In fact, in the case $p>0$, we have $\mathbb{E}\nu_N\sim Q_pN^{\frac{p}{p+1}}$ and it is not
  hard to see that, in agreement with the claim in \cite{mut}, $Q_pN^{\frac{p}{p+1}}\sim S(e^{-\sigma_N}),\
  N\rightarrow\infty$, where $\sigma_N$ is as in \refm[27].
  \\ A. Barbour and B. Granovsky (\cite{bar1}) explored the case when $Z_k$ in \refm[90] are
  quite general r.v.'s obeying $\mathbb{E}Z_k\sim k^{p-1},\ k\to \infty,\quad
  p<0.$ (This includes assemblies with
  $m_k\sim  k^{p-1}k!,\ k\to \infty,\ p<0$). It was shown in (\cite{bar1})
   that such RCS's
  exhibit a completely different asymptotic behavior. In particular, in
    this case $\nu_N$ is finite
  with probability $1$.
  \\ Hence, comparing the asymptotic behavior of  $\nu_N$ (e.g.  $\mathbb{E}\nu_N)$ in the cases $p<0$ and $p\geq0$, one sees
  that $p=0$ can be viewed as a point of phase transition for the measure $\mu_N$ as
  $N\to\infty$.
    \vskip .5cm
   \nin In conclusion, we provide a few examples of assemblies that conform to the setting of the present paper:
  \vskip .2cm
   {\bf{Example 1.$\;$Forests of labelled and colored linear
   trees.}}\quad A linear tree (see \cite{bur}) is a graph
   with no cycles, where each vertex has no more than two neighbors. Assuming that a vertex is labelled and
   is colored into one of $q\;(q\geq1)$ colors gives $m_k=q^kk!$,  $a_k=q^k,\ k\geq1$.  So, by
   the remark following Theorem 4.6, this RCS corresponds to the case $p=1$.

  \vskip .5cm
  {\bf{Example 2.$\;$Forests of labelled rooted linear trees.}}\quad In this case we have $m_k=kk!,\\ k\geq1$, which
   gives $a_k=k,\ k\geq1$, that corresponds to the case $p=2$.

 \vskip .5cm
  {\bf Example 3. Compositions.}$\;$(see \cite{sta})\quad  Consider a space $\Upsilon_N$ of ordered
  m-tuples$\\ {\xx}=(x_1,\ldots,x_m),\ m=1,\ldots,N$, where $x_i$  are positive integers, summing to $N$. In other
   words, $\Upsilon_N$ is a space of all ordered partitions (=compositions) of $N$. We define the probability measure
    $\lambda_N$ on $\Upsilon_N$:
  \be
   \la{333}
   \lambda_N({\xx})=\frac{1}{(m(\xx))!}{(r_N)}^{-1},\quad{\xx}=(x_1,\ldots,x_m)\in\Upsilon_N,
  \end{equation}
  where $({r_N})^{-1}$ is the normalizing constant and $m=m(\xx)$ is the number of components of $\xx$.
  This means that $\lambda_N$ prescribes
  to $\xx\in\Upsilon_N$ the weight $\frac{1}{(m(\xx))!}$. Denote $n_i=n_i({\xx}),\; i=1,\ldots,N$ the number of
  components of $\xx\in\Upsilon_N$ that are equal to $i$. Then we have
  \ber
   \la{999}
   r_N=\ds\sum_{m=1}^{N}{\frac{1}{m!}\sum_{\xx\in\Upsilon_N:\ m(\xx)=m}{1}}
   =\ds\sum_{m=1}^{N}{\frac{1}{m!}}\sum_{\eta\in\Omega_N:|\eta|=m}{\frac{m!}{n_1!\ldots n_m!}}=\ds\sum_{\eta\in\Omega_N}
   {\frac{1}{n_1!\ldots n_N!}},
 \ena
  where, as before, $\Omega_N$ is the set of all (unordered) partitions of $N$, while $|\eta|=n_1+\ldots+n_N$. \refm[999]
  says that $r_N=c_N$, where $c_N$ is the partition function of the measure $\mu_N$ in the case
   $a_k=1\:,\;k=1,\ldots,N$. Therefore, we can apply our results for $p=1$ to the number of summands in
   the random composition drawn according to $\lambda_N$.

\newpage
 \nin {\bf Acknowledgement}

\nin The research of the second author was supported by the Fund
for the Promotion of Research at Technion.

\nin The authors are thankful to Prof. Andrew Barbour for his
remarks on the first draft of the paper and to referees for their
stimulating and constructive suggestions.

\end{document}